\documentclass{amsart}
\usepackage[margin=1.5in]{geometry}

\usepackage[normalem]{ulem}
\usepackage{amsmath,amsthm, latexsym, amssymb,mathrsfs,  tikz-cd}
\usepackage{stackengine,scalerel, quiver, MnSymbol}
\usepackage[colorlinks]{hyperref}
\usepackage{cleveref}
\definecolor{dark-red}{rgb}{0.6,0,0}
\definecolor{dark-green}{rgb}{0,0.4,0}
\definecolor{medium-blue}{rgb}{0,0,0.5}
\hypersetup{
    colorlinks, linkcolor={dark-red},
    citecolor={dark-green}, urlcolor={medium-blue}
}

\usepackage{relsize}
\usepackage[bbgreekl]{mathbbol}
\usepackage{amsfonts}
\DeclareSymbolFontAlphabet{\mathbb}{AMSb} %to ensure that the meaning of \mathbb does not change
\DeclareSymbolFontAlphabet{\mathbbl}{bbold}

\usepackage{raysmacros}

\newcommand{\Gr}{\mathrm{Gr}}

\newcommand{\Spd}{\mathrm{Spd}}

\newcommand{\colim}{\mathrm{colim}}

\newcommand{\mc}[1]{\mathcal{#1}}

\newcommand{\mr}[1]{\mathrm{#1}}
\newcommand{\mf}[1]{\mathfrak{#1}}

\newcommand{\dR}{\mathrm{dR}}
\newcommand{\BC}{\mathrm{BC}}

\newcommand{\Spa}{\mathrm{Spa}}

\newcommand{\FF}{\mathrm{FF}}
\DeclareMathOperator{\Hom}{Hom}

\newcommand{\crys}{\mr{crys}}

\numberwithin{equation}{subsection}
\numberwithin{equation}{subsubsection}

\theoremstyle{plain}

\newtheorem{maintheorem}{Theorem} 
 
\newtheorem{maincorollary}[maintheorem]{Corollary}

\newtheorem*{theorem*}{Theorem}
\newtheorem*{conjecture*}{Conjecture}

\newtheorem{theorem}[subsubsection]{Theorem}

\newtheorem{proposition}[subsubsection]{Proposition}
\newtheorem{lemma}[subsubsection]{Lemma}

\theoremstyle{definition}

\newtheorem{example}[subsubsection]{Example}

\newtheorem{definition}[subsubsection]{Definition}
\newtheorem{remark}[subsubsection]{Remark}
\newtheorem{question}[subsubsection]{Question}

\newcommand{\BCfrac}[2]{\mathrm{BC}(\tfrac{#1}{#2})}
\newcommand{\dg}{\textnormal{\dag}}

\title[Intersections of analytic manifolds with Banach--Colmez spaces]{$\mathbb{Q}_p$-analytic intersections of $\mathbb{C}_p$-analytic manifolds with perfectoid Banach--Colmez spaces}

\author{Sean Howe and Ray Li}

\begin{document}

\begin{abstract} We show that transverse intersections between smooth $\mathbb{C}_p$-analytic manifolds in $\mathbb{C}_p^n$ and perfectoid Banach--Colmez spaces of complementary $\mathbb{C}_p$-dimension are $n$-dimensional $\mathbb{Q}_p$-analytic manifolds. In the case of the Banach--Colmez space $\BCfrac{1}{2}$, this is a strong form of the relativistic $p$-adic sunscreen conjecture. These results hold more generally with $\mathbb{C}_p$ replaced by an arbitrary complete algebraically closed extension of $\mathbb{Q}_p$, and we also establish some non-linear analogs for infinite level Lubin--Tate spaces.
\end{abstract}

\maketitle

\tableofcontents

\section{Introduction}

Let $C/\mathbb{Q}_p$ be a complete algebraically closed extension (for example, take $C=\mathbb{C}_p:=\overline{\mathbb{Q}}_p^\wedge$). There is an associated Fargues--Fontaine curve $\FF:=\FF_{C^\flat}$, an adic space over $\mathbb{Q}_p$, equipped with a canonical point $\infty: \Spa\, C \hookrightarrow \FF$. Let $\mathcal{E}$ be a vector bundle of rank $n$ and degree $d$ on $\FF$ with Harder--Narasimhan slopes lying in $[0,1)$, and let $\BC(\mathcal{E}) := H^0(\FF, \mathcal{E})$, a $\mathbb{Q}_p$-Banach space. We write $E:=H^0(\Spa\, C, \infty^*\mathcal{E})$, the fiber of $\mathcal{E}$ at $\infty$. Note that $E$ is an $n$-dimensional $C$-vector space and thus has a canonical topology; forgetting the $C$-linearity, it can be thought of also as a $\mathbb{Q}_p$-Banach space. Evaluation of global sections is a $\mathbb{Q}_p$-linear map $\iota: \BC(\mathcal{E}) \hookrightarrow E$ and identifies $\BC(\mathcal{E})$ with a closed subspace of $E$ (Lemma \ref{lem.HN-inj}).  

Let $M \subseteq E$ be a $C$-analytic submanifold --- for example, fixing an isomorphism $E=C^n$, $M$ could be $X(C)$ for $X \subseteq \mathbb{A}^n_C$ a smooth rigid analytic subvariety. We say $M$ intersects $\BC(\mathcal{E})$ transversely if, for every $m \in M \cap \BC(\mathcal{E})$, $T_{M,m} + \BC(\mathcal{E}) = E.$  

\begin{remark}
    If we interpret $\BC(\mathcal{E})$ as a $\mathbb{Q}_p$-analytic Banach manifold modeled on itself, then, for any $b  \in \BC(\mathcal{E})$, $T_{\BC(\mathcal{E}),b}=\BC(\mathcal{E})$. Thus this is the usual notion of a transverse intersection.  
\end{remark}

\begin{example} If $M$ has codimension $d$ (the degree of $\mathcal{E}$), then $T_{M,m} + \BC(\mathcal{E})=E$ if and only if $T_{M,m}$ lies in the $\mathcal{E}$-admissible locus of the Grassmannian $\Gr(n-d, E)$ of $(n-d)$-planes in $E$, i.e., the locus of subspaces $W \subseteq E$ such that modifying $\mathcal{E}$ by $W$ (i.e., forming the bundle whose sections are those sections of $\mathcal{E}$ that evaluate into $W$ at $\infty$) yields the trivial bundle --- see Lemma \ref{lemma.transv-adm}. If $\mathcal{E}$ is isoclinic of slope $1/n$ so that $\BC(\mathcal{E})=\BCfrac{1}{n}$, then $d=1$, the Grassmannian is $\mathbb{P}(E^*)$, and the admissible locus is all of $\mathbb{P}(E^*)$ so that the intersection is always transverse! In general, the admissible locus is open, so one can replace any $M$ of codimension $d$ with an open (but possibly empty) submanifold to obtain a transverse intersection with $\BC(\mathcal{E})$. 
\end{example}

When the codimension of $M$ is $d$ (the degree of $\mathcal{E}$), the same analysis in terms of admissibility shows that, for any $m \in M \cap \BC(\mathcal{E})$ where the intersection is transverse, $T_{M,m} \cap \BC(\mathcal{E})$ is an $n$-dimensional $\mathbb{Q}_p$-subspace. We thus expect that $M \cap \BC(\mathcal{E})$ has a natural structure of an $n$-dimensional $\mathbb{Q}_p$-analytic manifold around such a point, and that is our main result. 

\begin{maintheorem}\label{maintheorem.transverse} Let $\mathcal{E}$ be a vector bundle of rank $n$ and degree $d$ on $\FF$ with Harder--Narasimhan slopes lying in $[0,1)$, and suppose $M \subseteq E:=H^0(\Spa\, C, \infty^*\mathcal{E})\cong C^n$ is a $C$-analytic submanifold of codimension $d$. If $M$ intersects $\BC(\mathcal{E})$ transversely, then $M \cap \BC(\mathcal{E})$ is an $n$-dimensional $\mathbb{Q}_p$-analytic submanifold of the $\mathbb{Q}_p$-analytic Banach manifold $\BC(\mathcal{E})$.  
\end{maintheorem}

\begin{remark}
Since $\BC(\mathcal{E})$ is a closed $\mathbb{Q}_p$-subspace of the $\mathbb{Q}_p$-Banach space $E$ and, by \cite[Proposition 10.5]{Schneider.NonarchimedeanFunctionalAnalysis}, any closed $\mathbb{Q}_p$-subspace is complemented, $M \cap \BC(\mathcal{E})$ being a $\mathbb{Q}_p$-analytic submanifold of $\BC(\mathcal{E})$ is equivalent to $M \cap \BC(\mathcal{E})$ being a $\mathbb{Q}_p$-analytic submanifold of $E$. 
\end{remark}

\begin{remark}\label{remark.slope-1}If $\mathcal{E}$ has slopes in $[0,1]$, then the conclusion of Theorem \ref{maintheorem.transverse} still holds with the same proof, except that now $\iota:\BC(\mathcal{E}) \rightarrow E$ has kernel a potentially non-trivial finite-dimensional $\mathbb{Q}_p$-vector space and thus $M \cap \BC(\mathcal{E})$ should be written instead as $\iota^{-1}(M)$.\end{remark}

In the case that $\mathcal{E}$ is isoclinic of slope $1/n$, we can identify $\BC(\mathcal{E})=\BCfrac{1}{n}$ with $(B^+_\crys)^{\varphi^n=p}$, where $B^+_\crys$ is Fontaine's crystalline period ring associated to $C$ and $\varphi$ is its Frobenius. For $\theta: B^+_\crys \rightarrow C$ the Fontaine map, there is a choice of basis of $E$ so that $\iota$ becomes the map $b \mapsto \Theta(b):=(\theta(b), \theta(\varphi(b)), \ldots, \theta(\varphi^{n-1}(b)))$. Theorem \ref{maintheorem.transverse} has the following explicit consequence (note that the $n=1$ case is not technically contained in Theorem \ref{maintheorem.transverse} but it can be obtained by a separate trivial argument or by the extension of Remark \ref{remark.slope-1}):

\begin{maincorollary}\label{Cor.Germ}Let $F \in C[[x_1, \ldots, x_n]]$ have constant term $0$ and non-vanishing gradient at $0$, i.e.,
\[ F(x_1, \ldots, x_n)=0+ F_1  + \ldots \textrm{ where $F_1=\sum_{i=1}^n a_i x_i$ and $a_i \neq 0$ for some $i$},\]
and suppose $F$ converges in a neighborhood of zero. There is an open neighborhood $U$ of $0$ in $(B^+_\crys)^{\varphi^n=p}$ such that the vanishing set in $U$ of $b \mapsto F(\Theta(b))$ is homeomorphic to $\mathbb{Z}_p^n$. 
\end{maincorollary}

For $n=2$ and $C=\mathbb{C}_p$, Corollary \ref{Cor.Germ} establishes the relativistic $p$-adic sunscreen conjecture of \cite{Howe.TheRelativisticPAdicSunscreenConjecture}. The rest of Corollary \ref{Cor.Germ}/Theorem \ref{maintheorem.transverse} establishes the immediate generalizations of the conjecture described in \cite[\S1.4]{Howe.TheRelativisticPAdicSunscreenConjecture}. Even in the $n=2$ case for the simplest non-linear example, $F(x,y)=y^2-x$, it was unclear to us before this work whether there was even a single point of $\BCfrac{1}{2}$ in the vanishing locus of $F$ besides $(0,0)$ (see \cite[Example 1.2.1]{Howe.TheRelativisticPAdicSunscreenConjecture}).

\begin{remark}Theorem \ref{maintheorem.transverse} and its proof imply a stronger conclusion than stated in Corollary \ref{Cor.Germ}: the vanishing set in $U$ is given by an embedding of $\mathbb{Q}_p$-analytic Banach manifolds $\mathbb{Z}_p^n \hookrightarrow U$, where here $U$ is interpreted as an infinite-dimensional $\mathbb{Q}_p$-analytic Banach manifold. 
\end{remark}

\begin{remark} Let us emphasize the unusual nature of the geometry here: we are intersecting an $(n-d)$-dimensional $C$-analytic manifold inside of $C^n$ with a linear infinite-dimensional $\mathbb{Q}_p$-analytic Banach manifold (note that $C$ is also infinite-dimensional over $\mathbb{Q}_p$) to obtain an $n$-dimensional $\mathbb{Q}_p$-analytic manifold --- see \cite{Howe.TheRelativisticPAdicSunscreenConjecture} for a flavorful take on this exotic structure.
\end{remark}

\begin{remark}\label{remark.BC-conn-perfd-mot}What we are calling $\BC(\mathcal{E})$ is the $C$-points of an effective Banach--Colmez space over $\Spd\, C$ that is usually denoted the same way. Our conditions on $\mathcal{E}$ guarantee that this effective Banach--Colmez space is perfectoid and embeds in $C^{\mathrm{rk} \mathcal{E}}$, explaining the title of this paper. As indicated in \cite{Howe.TheRelativisticPAdicSunscreenConjecture}, a major motivation for treating this linear case is as a step towards understanding intersections of perfectoid spaces and rigid analytic varieties in more general settings --- in particular, one would like to understand the geometry of the transverse intersections involving period maps and local Shimura varieties that are studied in \cite[\S6]{HoweKlevdal.AdmissiblePairsAndpAdicHodgeStructuresIIIVariationAndUnlikelyIntersection}. In that work, the notion of transversality used here is replaced by an algebraic theory based on nilpotents (the inscribed $v$-sheaves of \cite{Howe.InscriptionTwistorsAndPAdicPeriods}), and in that generality it is not clear how the structure of the Banach--Colmez Tangent Space is reflected in the local punctual geometry. In the present work we make partial progress for infinite level Lubin--Tate spaces --- see \S\ref{ss.non-linear-intro}. 

 \end{remark}

 \begin{question}\label{q.extension}
 In the setting of Theorem \ref{maintheorem.transverse} with $\mathcal{E}=\mathcal{O}(\frac{1}{n})$, suppose $M=X(C)$ for $X/C$ a smooth codimension 1 rigid analytic subvariety of $\mathbb{A}^n_C$ (fixing an identification $E=C^n=\mathbb{A}^n(C)$). If $C'/C$ is a larger complete algebraically closed extension then, writing $M'=X(C')$ and $\BC(\mathcal{E})'$ for the associated Banach--Colmez space over $C'$, we find that $M \cap \BC(\mathcal{E})$ is an open $\mathbb{Q}_p$-analytic submanifold of $M' \cap \BC(\mathcal{E})'$. Can $M' \cap \BC(\mathcal{E})'$ be strictly larger than $M \cap \BC(\mathcal{E})$?
 \end{question}

\subsection{Proofs}

Our proof of Theorem \ref{maintheorem.transverse} is via the theory of $\mathbb{Q}_p$-analytic Banach manifolds. However, Theorem \ref{maintheorem.transverse} can also be proved using the forthcoming overconvergent representability results of \cite{LiOcdt}. Rather than giving that proof in full generality, we content ourselves with using it to give another deduction of Corollary \ref{Cor.Germ}, which already illustrates the key ideas. 

\subsubsection{Sketch of proof via Banach manifolds.} 
We may work locally around any point $m \in M \cap \BC(\mathcal{E})$. Transversality gives a surjection $\BC(\mathcal{E}) \twoheadrightarrow E / T_{M,m}$, and the connection between transversality and admissibility shows the kernel $Q$ is an $n$-dimensional $\mathbb{Q}_p$-vector space. The key is then to choose a splitting of this surjection of $\mathbb{Q}_p$-Banach spaces --- this exists by Lemma \ref{lem.splitting}. After choosing this splitting, we may then work in the category of $\mathbb{Q}_p$-analytic Banach manifolds. There is an implicit function theorem in this category, and applying it yields the  result. It is thus crucial that we chose the right moment to let go of the finiteness properties of Banach--Colmez spaces and instead treat them just as $\mathbb{Q}_p$-Banach spaces. 

We note that the theory of non-archimedean analytic Banach manifolds was exposited without proofs in \cite{Bourbaki.VarietesDifferentiellesEtAnalytique1-7}, but it does not seem to have been discussed much since then (there is a reasonable literature on the complex analytic theory, but in the non-archimedean case it seems only the continuous and smooth theories found traction!). The proofs of the basic results, though unwritten, are essentially the same as in the finite-dimensional case (or in the better-documented infinite-dimensional case over the complex numbers). 

\subsubsection{Sketch of proof via the representability theorem} In forthcoming work \cite{LiOcdt}, one of the authors (RL) establishes a representability theorem for solutions of Maurer--Cartan equations in certain topological DGAs related to $p$-adic analytic deformations. In \S\ref{s.representability} we state an instance of this result (Theorem \ref{thm:Representability}), and explain how it can be used to deduce Corollary \ref{Cor.Germ} by first establishing an algebraic version. The quadratic case fits naturally in this framework, while the general case is more naturally explained using $A_\infty$-algebras.  

\subsection{The non-linear case}\label{ss.non-linear-intro}
In \S\ref{s.non-linear} we explain how the Banach manifold perspective can also be applied to the study of the infinite level Lubin--Tate spaces and their period mappings, giving a non-linear analog of the study of intersections in Banach--Colmez spaces discussed above. In particular, we explain how this connects to the theory of inscribed $v$-sheaves and the unlikely intersection conjectures of \cite[\S6]{HoweKlevdal.AdmissiblePairsAndpAdicHodgeStructuresIIIVariationAndUnlikelyIntersection}. We hope to develop these ideas further in future work using a more refined notion of Banach Manifolds (see Remarks \ref{remark.hk-intersections} and \ref{remark.BanachManifolds} and note that here we have invoked Colmez Capitalization). 

\subsection{Outline}
In \S\ref{s.preliminaries} we discuss some preliminaries on Banach spaces and Banach--Colmez spaces. In particular, we recall the Banach topology on a Banach--Colmez space of the form $\BC(\mathcal{E})$ and describe the linear computation relating transversality and admissibility. 

In \S\ref{s.banach-man}, we exposit the proof of Theorem \ref{maintheorem.transverse} via $\mathbb{Q}_p$-analytic Banach manifolds. Since the theory  of $p$-adic analytic Banach manifolds may be unfamiliar to readers, in \S\ref{s.an-example} we first illustrate the ideas behind Theorem \ref{maintheorem.transverse} by constructing non-trivial points in $\BCfrac{1}{2} \cap V(y^2-x)$ using the Banach fixed point theorem. In \S\ref{ss.ban-man} we recall the general theory of $p$-adic analytic Banach manifolds from \cite{Bourbaki.VarietesDifferentiellesEtAnalytique1-7} before giving the proof of Theorem \ref{maintheorem.transverse} in \S\ref{s.main-proof}.

In \S\ref{s.representability}, we exposit the proof of Corollary \ref{Cor.Germ} via the representability theorem of \cite{LiOcdt}. After stating the relevant version in \S\ref{ss.def-thy-prelim}, we first establish an algebraic analog in \S\ref{ss.alg-inc}, then deduce the topological statement in \S\ref{ss.top-version}. The approach via DGAs is most natural in the quadratic case and we explain this at the start of \S\ref{ss.alg-inc}. For general power series it is clearer to use $A_\infty$-algebras; this is thus the approach we emphasize in the remainder of \S\ref{ss.alg-inc}. However, since the results of \cite{LiOcdt} are given only in terms of DGAs, we supplement this discussion in \S\ref{appendix:DGAConstruction} by explicitly constructing  an equivalent DGA that can be used in the general case.  

In \S\ref{s.non-linear} we give our results on the non-linear case of the infinite level Lubin--Tate space. 

\subsection{Acknowledgements}SH was partially supported during the preparation of this work by the NSF under grants DMS-2201112, DMS-2501816, and DMS-2540284. 

The construction of points in $\BCfrac{1}{2} \cap V(y^2-x)$ presented in \S\ref{s.an-example} was developed in dialog with ChatGPT 5.5. Claude and ChatGPT were both used to proofread this paper.

\section{Preliminaries}\label{s.preliminaries}

\subsection{Banach spaces}\label{ss.ban-spa-prelim}
Let $L/\mathbb{Q}_p$ be an extension that is complete for an absolute value extending the $p$-adic absolute value on $\mathbb{Q}_p$. An $L$-Banach space is a complete locally convex topological $L$-vector space $B$ whose topology can be defined by a non-archimedean norm --- we note that the norm is not a part of the data, and that any two norms defining the same topology are equivalent. We will need the following result:
\begin{lemma}\label{lem.splitting}Suppose $L$ is spherically complete. If $f: V \twoheadrightarrow W$ is a continuous $L$-linear surjection between $L$-Banach spaces and $\ker f$ is a finite-dimensional $L$-vector space, then there exists a continuous $L$-linear section of $f$, $s: W \rightarrow V$.  
\end{lemma}
\begin{proof}
Because $L$ is spherically complete, the Hahn--Banach theorem holds \cite[Proposition 9.2]{Schneider.NonarchimedeanFunctionalAnalysis}. Thus, choosing a basis $e_1, \ldots, e_n$ for $\ker f$, we may extend the elements of the dual basis $\delta_1, \ldots, \delta_n$ of $(\ker f)^*$ to continuous linear functionals $\tilde{\delta}_i: V \rightarrow L$. The map $g=(\tilde{\delta}_1, \ldots, \tilde{\delta}_n): V \rightarrow L^n$ is then a left-splitting of $L^n=\ker f \hookrightarrow V$, so that $V=\ker f \oplus \ker g$. The restriction of $f$ to $\ker g$ is then a continuous $L$-linear bijection so, by the open mapping theorem, is an isomorphism of $L$-Banach spaces. Taking the inverse yields the section $s$.     
\end{proof}

\begin{remark}
In the proof of Lemma \ref{lem.splitting}, if $L$ is discretely valued, then one can replace the explicit use of the Hahn--Banach theorem above %to construct a complement
with the more general statement that any closed subspace is complemented \cite[Proposition 10.5]{Schneider.NonarchimedeanFunctionalAnalysis}.  
\end{remark}

\subsubsection{}\label{sss.analytic-maps}
In our discussion of $p$-adic analytic Banach manifolds we will use analytic functions on an open in one $p$-adic Banach space valued in another $p$-adic Banach space; we explain this now. Suppose $B_1$ and $B_2$ are two $L$-Banach spaces. For $k \geq 1$, a \emph{homogeneous polynomial function of degree $k$} from $B_1$ to $B_2$ is a function $f: B_1 \rightarrow B_2$ such that there exists a continuous symmetric multilinear map $m: B_1^k \rightarrow B_2$ such that $f(b_1)=m(b_1, \ldots, b_1)$ for all $b_1 \in B_1$. A homogeneous polynomial function of degree $0$ is a constant function $f: B_1 \rightarrow B_2$. Fixing norms defining the topologies on both $B_1$ and $B_2$ to be denoted as $|| \cdot ||$, for $f$ a homogeneous polynomial function, we write $||f||=\mathrm{sup}_{b \neq 0} ||f(b)||/||b||^k$. 

For $U \subseteq B_1$ an open subset, a function $f: U \rightarrow B_2$ is \emph{analytic at $u$} if there exist homogeneous polynomial functions $f_k: B_1 \rightarrow B_2$ of degree $k$, $k \geq 0$, and a real number $r >0$ such that $B_r(u)=\{v \in B_1\, |\, ||v-u||\leq r\} \subseteq U$, $\lim_{k \rightarrow \infty} ||f_k|| r^k=0$, and, for any $v \in B_r(u)$, 
\[ f(v) = \sum_{k \geq 0} f_k(v-u) \]
where the sum on the right is convergent by the limit condition (as the norm is ultrametric). Note that, by scaling $r$, this property is independent of the choice of norms.

If $f$ is analytic at $u$ then the derivative of $f$ exists at $u$ and is equal to $f_1$; in particular, $f_1$ is uniquely determined by $f$. We say $f$ is \emph{analytic} if it is analytic at every $u \in U$.

\subsection{Banach--Colmez spaces}
Let $C/\mathbb{Q}_p$ be a complete algebraically closed extension.

\subsubsection{}Attached to 
$C$ there is a Fargues--Fontaine curve $\FF:=\FF_{C^\flat}$, an adic space over $\mathbb{Q}_p$, and it is equipped with a canonical point $\infty: \Spa\,C \hookrightarrow \FF$. 

\subsubsection{}Recall that there is a slope formalism for vector bundles on $\FF$; moreover, for each $\lambda \in \mathbb{Q}$ there is a unique simple bundle $\mathcal{O}(\lambda)$ of slope $\lambda$ and any vector bundle is a direct sum of these simple bundles. For $\mathcal{E}$ a vector bundle, the Harder--Narasimhan slopes of $\mathcal{E}$ are precisely those $\lambda$ such that $\mathcal{O}(\lambda)$ appears as a summand of $\mathcal{E}$. We further recall that:
\begin{enumerate}
    \item $H^0(\FF, \mathcal{O})=\mathbb{Q}_p$, 
    \item for $\lambda \geq 0$, $H^1(\FF,\mathcal{O}(\lambda))=0$, and
    \item for $\lambda < 0$, $H^0(\FF, \mathcal{O}(\lambda))=0$ and $H^1(\FF, \mathcal{O}(\lambda)) \neq 0$.
\end{enumerate}

\subsubsection{}\label{sss.bc-banach} Let $\mathcal{E}$ be a rank $n$ vector bundle on $\FF$. We write $\BC(\mathcal{E}):=H^0(\FF, \mathcal{E})$. This is a $\mathbb{Q}_p$-vector space, but we claim it has a natural $\mathbb{Q}_p$-Banach space structure. Indeed, since $\FF$ is quasi-compact, we may choose a finite cover $U_i=\Spa(A_i, A_i^+)$, $i=1, \ldots, k$, and trivializations $\mathcal{E}|_{U_i} \xrightarrow{\sim} \mathcal{O}_{U_i}^n$. Then we may take as a unit ball in $\BC(\mathcal{E})$ those sections whose restriction to each $U_i$ lies in $(A_i^+)^n \subseteq A_i^n = H^0(U_i, \mathcal{E})$ (which is a unit ball in $A_i^n$ since $\FF$ is uniform). Different choices give rise to the same $\mathbb{Q}_p$-Banach topology on $\BC(\mathcal{E})$. 

\begin{remark}Regarding $\BC(\mathcal{E})$ as a $\mathbb{Q}_p$-Banach space is much weaker than regarding it as a Banach--Colmez space. Indeed, viewing objects in the latter category as functors on perfectoid spaces over $C$, this amounts to restricting such a functor to profinite sets (viewed as perfectoid spaces over $C$ by sending $T$ to $\Spa(\mathrm{Cont}(T,C), \mathrm{Cont}(T, \mathcal{O}_C))$). 
\end{remark}

\subsubsection{} There is a natural $\mathbb{Q}_p$-linear restriction map $\iota: \BC(\mathcal{E})\rightarrow H^0(\Spa\, C, \infty^*\mathcal{E}) \cong C^n$. It is evidently continuous for the canonical topology on the finite-dimensional $C$-vector space $H^0(\Spa\, C, \infty^*\mathcal{E})$ and the topology described above on $\BC(\mathcal{E})$. We also have
\begin{lemma}\label{lem.HN-inj}
If the Harder--Narasimhan slopes of $\mathcal{E}$ are $< 1$, then $\iota$ is a closed embedding, i.e., it is injective and identifies $\BC(\mathcal{E})$ with a closed $\mathbb{Q}_p$-subspace of $H^0(\Spa\, C, \infty^*\mathcal{E})$.  
\end{lemma}
\begin{proof}
    For the injectivity, by passing to direct summands, it suffices to treat the case $\mathcal{E}=\mathcal{O}(\lambda)$. As the ideal sheaf of $\infty$ is isomorphic to $\mathcal{O}(-1)$, the kernel of $\iota$ is 
    \[ H^0(\FF, \mathcal{O}(\lambda)\otimes \mathcal{O}(-1))=H^0(\FF, \mathcal{O}(\lambda -1))=0\]
    where the last equality is because $\lambda -1 <0$. 

    We now show $\iota$ identifies $\BC(\mathcal{E})$ with a closed $\mathbb{Q}_p$-subspace. Since $\iota$ is continuous and $\mathbb{Q}_p$-linear, the open mapping theorem for $\mathbb{Q}_p$-Banach spaces implies it suffices to show that the image of $\iota$ is closed. But, by the long exact sequence in cohomology attached to the short exact sequence
    \[ 0 \rightarrow \mathcal{E}(-\infty) \rightarrow \mathcal{E} \rightarrow \infty_*\infty^*\mathcal{E} \rightarrow 0, \]
    this image is the kernel of the connecting map 
    \[ \delta: H^0(\Spa\, C, \infty^*\mathcal{E}) \rightarrow H^1(\FF, \mathcal{E}(-\infty)). \]
    We claim the codomain has a natural Hausdorff topology for which $\delta$ is continuous, implying that the kernel is a closed subspace. One way to obtain this is by viewing both the domain and codomain as condensed sets (by restricting the evident extension to all perfectoid spaces over $C$ to profinite sets over $C$), and noting that the underlying topological space of the first is the finite-dimensional $C$-vector space $H^0(\Spa\, C, \infty^*\mathcal{E})$ with its usual topology while the underlying topology of the second is the Banach--Colmez topology on $H^1(\FF, \mathcal{E}(-\infty))$. 
\end{proof}

\subsubsection{}
For $W\subseteq H^0(\Spa\, C, \infty^*\mathcal{E})$ a $C$-subspace, we may form the modified vector bundle $\mathcal{E}_W$ that agrees with $\mathcal{E}$ away from $\infty$ and whose sections over a neighborhood of $\infty$ are those sections of $\mathcal{E}$ whose evaluation at $\infty$ lies in $W$. We say $W$ is \emph{$\mathcal{E}$-admissible} if $\mathcal{E}_W$ is a trivial bundle, i.e., isomorphic to $\mathcal{O}_{\FF}^n$ (equivalently, semistable of slope 0).

\begin{lemma}\label{lemma.adm-fd}If $W$ is $\mathcal{E}$-admissible, then the pre-image of $W$ under the restriction map $\BC(\mathcal{E}) \rightarrow H^0(\Spa\, C, \infty^*\mathcal{E})$ is an $n$-dimensional $\mathbb{Q}_p$-vector space. 
\end{lemma}
\begin{proof}
    Indeed, the pre-image of $W$ is precisely $\BC(\mathcal{E}_W)$, and $H^0(\FF, \mathcal{O}_\FF^n)=\mathbb{Q}_p^n$.  
\end{proof}

The following lemma establishes the link between transversality and admissibility.

\begin{lemma}\label{lemma.transv-adm}Suppose the Harder--Narasimhan slopes of $\mathcal{E}$ are in $[0,1)$ and write $E:=H^0(\Spa\, C, \infty^*\mathcal{E})$. Viewing $\BC(\mathcal{E})$ as a subspace of $E$ via Lemma \ref{lem.HN-inj}, a $C$-subspace $W\subseteq E$ of codimension equal to the degree of $\mathcal{E}$ is $\mathcal{E}$-admissible if and only if $W + \BC(\mathcal{E}) = E.$
\end{lemma}
\begin{proof}
Recalling $E:=H^0(\Spa\, C, \infty^*\mathcal{E})$, the modification construction gives an exact sequence 
\[ 0 \rightarrow \mathcal{E}_W \rightarrow \mathcal{E} \rightarrow \infty_* E/W \rightarrow 0 \]
where the last map is evaluation of sections at $\infty$ composed with projection $E \rightarrow E/W$. Taking global sections, we obtain an exact sequence
\[ 0 \rightarrow \BC(\mathcal{E}_W) \rightarrow \BC(\mathcal{E}) \rightarrow E/W \rightarrow H^1(\FF, \mathcal{E}_W) \rightarrow 0 \]
where we used that $H^1(\FF, \mathcal{E})=0$ because $\mathcal{E}$ has Harder--Narasimhan slopes $\geq 0$. Thus $W + \BC(\mathcal{E})=E$ if and only if $H^1(\FF, \mathcal{E}_W)=0$. Since $\mathcal{E}_W$ is of slope $0$ (its degree is the degree of $\mathcal{E}$ minus the codimension of $W$), the latter happens if and only if $\mathcal{E}_W$ is a trivial bundle (indeed, otherwise for it to be of slope zero it must contain summands of the form $\mathcal{O}(a)$ and $\mathcal{O}(b)$ with $a<0$ and $b>0$, but then $H^1(\FF, \mathcal{O}(a))\neq 0$). 
\end{proof}

\begin{example}\label{example.bc1n}
    If $\mathcal{E}=\mathcal{O}(\frac{1}{n})$, $n > 1$, then \emph{any} codimension $1$ subspace of $E:=H^0(\Spa\, C, \infty^*\mathcal{E})$ is $\mathcal{E}$-admissible (this is equivalent to the surjectivity of the Gross--Hopkins/Hodge period map). It then follows from Lemma \ref{lemma.transv-adm} and Lemma \ref{lemma.adm-fd} that, for \emph{any} non-zero linear functional $f: E \rightarrow C$, $f\circ \iota: \BC(\mathcal{E}) \rightarrow C$ is surjective with kernel an $n$-dimensional $\mathbb{Q}_p$-vector space. For $n=2$, this is equivalent to the $p$-adic sunscreen property described in \cite{Howe.TheRelativisticPAdicSunscreenConjecture}. 
\end{example}

\section{Proof via $\mathbb{Q}_p$-analytic Banach manifolds}\label{s.banach-man}

In \S\ref{s.an-example} we explain how to use the Banach fixed point theorem to construct points in $\BCfrac{1}{2} \cap V(y^2-x)$. This is essentially a special case of the proof of the implicit function theorem, whose statement in the case of $p$-adic analytic Banach manifolds is recalled in \S\ref{ss.ban-man}.  In \S\ref{s.main-proof} we then use this implicit function theorem to prove Theorem \ref{maintheorem.transverse}. 

\subsection{An example}\label{s.an-example}

We fix a complete algebraically closed extension $C/\mathbb{Q}_p$ and we view the Banach--Colmez space $\BCfrac{1}{2}$ as a $\mathbb{Q}_p$-Banach space. It is equipped with a natural continuous $\mathbb{Q}_p$-linear injection
$\iota=(\iota_1, \iota_2):  \BCfrac{1}{2} \hookrightarrow C^2$ (see Lemma \ref{lem.HN-inj}). In this subsection we construct an injection from $\mathbb{Z}_p^2$ into $\BCfrac{1}{2} \cap V(F)$ for $F(x,y)=y^2-x$. 

\subsubsection{The $p$-adic sunscreen property}
By Example \ref{example.bc1n}, for any non-zero $C$-linear functional $f: C^2 \rightarrow C$, $f \circ \iota$ is surjective with kernel a 2-dimensional $\mathbb{Q}_p$-vector space. Applying this to $f$ the projection onto the first coordinate, $f(x,y)=x$, we find $\iota_1$ is a surjection $\BCfrac{1}{2} \twoheadrightarrow C$ and, writing $Q=\ker \iota_1$, $Q$ is a 2-dimensional $\mathbb{Q}_p$-subspace of $\BCfrac{1}{2}$. 

\subsubsection{Construction}\label{ss.construction}
We treat $\BCfrac{1}{2}$ and $C$ as $\mathbb{Q}_p$-Banach spaces. Via Lemma \ref{lem.splitting}, we choose a continuous $\mathbb{Q}_p$-linear section $s: C \rightarrow \BCfrac{1}{2}$ of $\iota_1$. By construction, $\iota_1( s(\alpha))=\alpha$ for any $\alpha \in C$, but we emphasize that $\iota_2 \circ s$ may not be $C$-linear (but it is $\mathbb{Q}_p$-linear).

The section $s$ induces a decomposition $\BCfrac{1}{2}=Q \oplus s(C)$. Since $s$ is continuous, we can and do choose an $A \geq 1$ such that $|\iota_2(s(\alpha))|\leq A|\alpha|$ for all $\alpha \in C$. In this decomposition, 
\[ F(\iota(q+s(\alpha)))= \iota_2(q + s(\alpha))^2 - \iota_1(q+s(\alpha))=\iota_2(q + s(\alpha))^2 - \alpha.\]

Setting  
\[ \Phi_q(\alpha):=\iota_2(q+s(\alpha))^2, \] 
we find $F(\iota(q+s(\alpha)))=0$ if and only if $\alpha$ is a fixed point of $\Phi_q$. Now, we claim that for $\epsilon = \frac{1}{2A^2}$, if $|\iota_2(q)| \leq \epsilon$ then $\Phi_q$ is a contraction on $B_\epsilon:=\{\alpha \in C\,|\, |\alpha| \leq \epsilon\}$. Indeed, if $|\iota_2(q)|\leq \epsilon$, then 
\begin{align*} |\Phi_q(\alpha)-\Phi_q(\beta)| &= |\iota_2(q+s(\alpha))^2 - \iota_2(q+s(\beta))^2| \\
 &=|2\iota_2(q)(\iota_2(s(\alpha))-\iota_2(s(\beta))) + \iota_2(s(\alpha))^2 -\iota_2(s(\beta))^2| \\
&=|2\iota_2(q)\iota_2(s(\alpha-\beta)) + \iota_2(s(\alpha-\beta))\iota_2(s(\alpha+\beta))| \\
&\leq \max(\epsilon A|\alpha-\beta|,  A^2\epsilon|\alpha-\beta|) \leq \frac{1}{2}|\alpha - \beta|. \end{align*}
In particular the map $\Phi_q$ is contracting. This also implies $\Phi_q$ preserves $B_\epsilon$ as,  for $\alpha \in B_\epsilon$, we find $|\Phi_q(\alpha)-\Phi_q(0)| \leq \frac{1}{2}|\alpha|<\epsilon$, but $\Phi_q(0)=\iota_2(q)^2$ satisfies $|\Phi_q(0)|\leq \epsilon^2 < \epsilon$ (as $\epsilon < 1$), so $\Phi_q(0)$ is also a center of $B_\epsilon$. Thus, for any $q$ with $|\iota_2(q)|\leq\epsilon$,  the Banach fixed point theorem implies $\Phi_q$ has a unique fixed point $\alpha_q \in B_\epsilon$. 

In particular, if we fix an isomorphism $Q\cong \mathbb{Q}_p^2$ such that $|\iota_2(q)|\leq \epsilon$ for $q \in \mathbb{Z}_p^2$ (which is possible by the continuity of $\iota_2$), then we obtain our injection from $\mathbb{Z}_p^2$ into $\BCfrac{1}{2} \cap V(F)$ by sending $q \in \mathbb{Z}_p^2$ to $q+ s(\alpha_q)$.

\subsection{$L$-analytic Banach manifolds}\label{ss.ban-man}
Let $L/\mathbb{Q}_p$ be an extension that is complete for an absolute value extending the $p$-adic absolute value on $\mathbb{Q}_p$. 

\subsubsection{}An \emph{$L$-analytic Banach manifold} is a topological space $M$ equipped with an atlas of charts $(U_M, U_B, B, \varphi)$ where $U_M \subseteq M$ is an open set, $B$ is an $L$-Banach space, $U_B \subseteq B$ is an open set, and $\varphi: U_M \rightarrow U_B$ is a homeomorphism, subject to the constraint that all transition maps are analytic in the sense of \S\ref{sss.analytic-maps}. 

\subsubsection{}An $L$-analytic Banach manifold is finite-dimensional if, for every chart $(U_M, U_B, B, \varphi)$, $B$ is a finite-dimensional $L$-vector space (whose topology is then uniquely determined). The theory of finite-dimensional $L$-analytic Banach manifolds is the usual theory of $L$-analytic manifolds as treated, e.g., in \cite[Part II, Chapter III]{Serre.LieAlgebrasAndLieGroups};  for this reason we will refer to finite-dimensional $L$-analytic Banach manifolds simply as $L$-analytic manifolds.

\begin{example}
    Let $K/L$ be a complete extension. Since any $K$-Banach space is also an $L$-Banach space and any $K$-analytic function in the sense of \S\ref{sss.analytic-maps} is also $L$-analytic, there is a faithful forgetful functor from the category of $K$-analytic Banach manifolds to $L$-analytic Banach manifolds. In particular, when $K/L$ is an infinite extension, (finite-dimensional) $K$-analytic manifolds are also examples of infinite-dimensional $L$-analytic Banach manifolds. 
\end{example}

\subsubsection{}\label{sss.submanifold} If $M$ is an $L$-analytic Banach manifold, then a locally closed subspace $N \subseteq M$ is a \emph{submanifold} if, for every $x \in N$, there exists a chart $(U_M, U_B, B, \varphi)$ of $M$ with $x \in U_M$ and a decomposition $B=B_1 \oplus B_2$ into complementary closed subspaces such that 
\[ \varphi(N \cap U_M)=U_B \cap B_1. \]
Any submanifold of $M$ admits a natural $L$-analytic Banach manifold structure whose charts are obtained by taking, for any $(U_M, U_B, B, \varphi)$ as above, $(N \cap U_M, U_B \cap B_1, B_1, \varphi|_{N \cap U_M}).$ 

\subsubsection{}If $M$ is an $L$-analytic Banach manifold, then, for any point $m \in M$, there is a tangent space $T_{M,m}$ which is an $L$-Banach space. It can be identified with $B$ for any chart $(U_M, U_B, B, \varphi)$ such that $m \in U_M$ --- two different choices are identified by applying the derivative of the corresponding analytic transition map (see \S\ref{sss.analytic-maps}). 

\subsubsection{} If $M_1$ and $M_2$ are $L$-analytic Banach manifolds, then there is a product $L$-analytic Banach manifold $M_1 \times M_2$ whose underlying topological space is the product and whose charts are given by the products of charts on $M_1$ and $M_2$. In particular, for $m=(m_1, m_2) \in M_1 \times M_2$, $T_{M_1 \times M_2,m}=T_{M_1, m_1} \oplus T_{M_2, m_2}.$  

\subsubsection{}The following is the implicit function theorem of \cite[Th\'{e}or\`{e}me 5.6.7 on p. 45]{Bourbaki.VarietesDifferentiellesEtAnalytique1-7}:
\begin{theorem}\label{thm.implicit-function} Let $X_1, X_2$ and $X'$ be $L$-analytic Banach manifolds. Suppose $U \subseteq X_1 \times X_2$ is open, $f:U \rightarrow X'$ is analytic, and suppose $u=(x_1, x_2) \in U$ is such that, writing $T_{U,u}=T_{X_1, x_1} \oplus T_{X_2, x_2}$, $df_u|_{T_{X_2, x_2}}: T_{X_2, x_2} \rightarrow T_{X', f(u)}$ is an isomorphism. Then, for $Z:=\{v \in U\, |\, f(v)=f(u)\}$, there exist opens $V_1 \subseteq X_1$ and $V_2 \subseteq X_2$ such that $u \in V_1 \times V_2 \subseteq U$ and $Z \cap V_1 \times V_2$ is the graph of an analytic function $g: V_1 \rightarrow V_2$.   
\end{theorem}

\subsection{Proof of the main result}\label{s.main-proof}

\begin{proof}[Proof of Theorem \ref{maintheorem.transverse}]
    As in the statement of the theorem, we let $\mathcal{E}$ be a vector bundle of rank $n$ and degree $d$ and Harder--Narasimhan slopes in $[0,1)$. We view $\BC(\mathcal{E})$ as a $\mathbb{Q}_p$-Banach space (and thus also a $\mathbb{Q}_p$-analytic Banach manifold). We fix an isomorphism $H^0(\Spa\, C, \infty^*\mathcal{E})=C^n$ so that $\BC(\mathcal{E}) \subseteq C^n$. We let $M$ be a $C$-analytic submanifold of $C^n$ of dimension $n-d$ such that, for every $m \in M \cap \BC(\mathcal{E})$, $T_{M,m} + \BC(\mathcal{E})=C^n$. 
    
    We want to show that $M \cap \BC(\mathcal{E})$ is an $n$-dimensional $\mathbb{Q}_p$-analytic submanifold of $\BC(\mathcal{E})$ in the sense of \S\ref{sss.submanifold}. We will accomplish this by exhibiting, for each $m \in M \cap \BC(\mathcal{E})$, a neighborhood of $m$ in $\BC(\mathcal{E})$ that is isomorphic to $V \times W$ where $V$ is open in $\mathbb{Q}_p^n$ and $W$ is open in $C^d$ such that $M \cap V \times W = V \times \{w\}$ for a point $w \in W$.  

    To that end, fix $m \in M \cap \BC(\mathcal{E})$. There is a neighborhood $U$ of $m$ in $C^n$ such that $U \cap M$ is the vanishing locus of a $C$-analytic submersion $f: U \rightarrow C^d$. In particular, writing $U':=U \cap \BC(\mathcal{E})$, $M \cap U'$ is the vanishing locus of $f|_{U'}$, a $\mathbb{Q}_p$-analytic function $U' \rightarrow C^d$. By our assumptions, Lemma \ref{lemma.transv-adm} and Lemma \ref{lemma.adm-fd} show  that $df_m|_{\BC(\mathcal{E})}$ is surjective with kernel an $n$-dimensional $\mathbb{Q}_p$-vector space $Q$. Thus, using Lemma \ref{lem.splitting}, we can and do choose a continuous $\mathbb{Q}_p$-linear section $s: C^d \rightarrow \BC(\mathcal{E})$ to write $\BC(\mathcal{E})=Q \oplus s(C^d)= Q \times s(C^d)$.

    Then $d(f|_{U'})_m|_{s(C^d)}$ is the isomorphism $s(C^d) \rightarrow C^d$, so Theorem \ref{thm.implicit-function} applies, and we find there are opens $V_1 \subseteq Q$ and $V_2 \subseteq s(C^d)$ such that $m \in V_1 \times V_2 \subseteq U' \subseteq \BC(\mathcal{E})$ 
    and a $\mathbb{Q}_p$-analytic function $g: V_1 \rightarrow V_2$ such that $M \cap V_1 \times V_2$ is the graph of $g$. In particular, after passing to an open $V \subseteq V_1$ and taking $W$ a small ball around 0 in $C^d$, $\psi: (q, z) \mapsto (q, g(q) + s(z))$ is an isomorphism of $V \times W$ onto a neighborhood of $m$ in $\BC(\mathcal{E})$ such that $\psi^{-1}(M)=V \times \{0\}$, thus we conclude. 

\end{proof}

\section{Proof via the representability theorem}\label{s.representability}

In this section we use the representability results of \cite{LiOcdt} to prove Corollary \ref{Cor.Germ}; we also prove an algebraic version in families,
see \Cref{prop:AlgebraicPadicSunscreen} for a precise statement.

We now sketch the proof since it borrows methods from an area
that may not be well known to the likely audience.
Our strategy uses the tools of homological perturbation theory/homotopy transfer.
To apply these, we view $F$ as a perturbation of a two-term chain complex.

First, suppose $F$ is linear.
Then, the desired set is $H^1$ of the complex concentrated in degrees $[1, 2]$
$$
\left[
    \BCfrac{1}{n}
    \xrightarrow{F \circ \Theta} C
\right].
$$

For more general $F$ we perturb this complex.
The next easiest case is when $F$ is quadratic.
Split $F = F_1 + F_2$ into its linear and quadratic part.
Consider
$$
\mb{T}
:=
\left[
    \BCfrac{1}{n}
    \xrightarrow{F_1 \circ \Theta} C
\right],
$$
as a chain complex concentrated in degrees $[1,2]$.
Then one can, in a straightforward manner, endow $\mb{T}$ with a differential graded algebra (DGA) structure
such that
$$
\mrm{MC}(\mb{T})
:=
\{
	x \in \BCfrac{1}{n}
	\mid
	dx - x \cdot x = 0
\}
$$
bijects with
$\{x \in \BCfrac{1}{n} \mid F(\Theta(x)) = 0\}$.
For example, when $n = 2$ and $F = x - y^2$,
$x \cdot x := \theta(\varphi(x)) \cdot \theta(\varphi(x))$.
The set $\mrm{MC}(\mb{T})$ is the set of solutions to the \emph{Maurer--Cartan equation} of $\mb{T}$
and is a fundamental object in deformation theory.

More generally, the degree $m$ part of $F$ should be encoded as an $m$-ary operation.
This endows $\mb{T}$ with the structure of a homotopy associative algebra, also known as an $A_\infty$-algebra.
The Maurer--Cartan equation can be generalized to this setting
and its solutions again recover $\{x \in \BCfrac{1}{n} \mid F(\Theta(x)) = 0\}$.

Hence we have reduced the problem to the study of the structure of $\mrm{MC}(\mb{T})$.
This set has a key property: it is a homotopy invariant in a neighborhood of $0$, for a suitable notion of homotopy
(here we use that $\mb{T}^0 = 0$;
more generally the gauge equivalence classes form a homotopy invariant).
The theory of homotopy transfer/homological perturbation
states that the cohomology of $\mb{T}$, $H^{\bullet}$,
can be endowed with an $A_\infty$-algebra structure such that $\mb{T} \sim H^{\bullet}$
for this suitable notion.
But since $H^{\bullet}$ is concentrated in a single degree,
$\mrm{MC}(H^{\bullet}) \cong H^1 \cong \mb{Q}_p^n$
from which Corollary \ref{Cor.Germ} follows.

% Maybe cite: https://arxiv.org/pdf/0710.2645
% https://arxiv.org/pdf/math/9906036 ?

\subsection{Deformation theory preliminaries}\label{ss.def-thy-prelim}
To make this precise,
we use the main representability theorem, Theorem 1.4, of \cite{LiOcdt},
specialized to the Banach case (\cite[Remark 2.28]{LiOcdt}):

First we need to set some notation.

\begin{definition}
    Let $\mb{Q}_p\{X_1,\ldots,X_m\} := \colim_{k \to \infty} \mb{Q}_p\langle X_1/p^k, \ldots, X_m/p^k\rangle$.
\end{definition}

These rings are Noetherian \cite[4.1 Proposition 6]{BoschLectures}.
Hence the $I$ in the following definition
are finitely generated and are defined
over some term in the colimit defining
$\mb{Q}_p\{X_1,\ldots,X_m\}$.

\begin{definition}
    Let $\mc{C}^{\dagger}$
    denote the category of analytic/overconvergent germs of affinoid algebras,
    i.e., rings of the form
    $
    R = \mb{Q}_p\{T_1, \ldots, T_m\}/I,
    $
    where $I \trianglelefteq \mb{Q}_p\{T_1, \ldots, T_m\}$ is a proper ideal.
    Morphisms in $\mc{C}^\dagger$
    are $\mathbb{Q}_p$-algebra morphisms.
    
    Every such $R$ is local.
    Let $\mathrm{Nil}^\dagger(R) := \langle T_1, \ldots, T_m\rangle$
    denote the maximal ideal.
\end{definition}

All tensors henceforth are solid tensors over $\mb{Z}_p$.

\begin{definition}
	Let $\mb{T}$ be a $\mb{Q}_p^{\blacksquare}$ differential graded algebra.
	Let $(R, f:R \to \mb{Q}_p)$ be an augmented $\mb{Q}_p^{\blacksquare}$-algebra.
    Let $\mf{m}_R := \ker(f)$.
	The \emph{set of Maurer--Cartan solutions} in $R$ is
	$$
	\mathrm{MC}(\mb{T},R) := \{x \in \mb{T}^1 \otimes \mf{m}_R \mid dx - x \cdot x = 0\},
	$$
    where $dx - x \cdot x \in \mathbb{T}^2 \otimes \mf{m}_R$.
\end{definition}

In all our usages, $\mf{m}_R$ will be the ideal generated by $T_1, \ldots, T_m$.
Hence, we have slightly abusively suppressed
$f$ from the notation.

The version of the representability theorem
we will use is the following.

\begin{theorem}
	\label[theorem]{thm:Representability}
	Let $\mb{T}$ be a differential graded algebra over $\mb{Q}_p^{\blacksquare}$.
	Suppose that $\mb{T}^i$ is a Banach space for each $i$
	(and $d$ and the multiplication are continuous).
	Assume that the cohomology of $\mb{T}$, $H^\bullet(\mb{T})$,
	is finite-dimensional in each degree.
	Further suppose that $\mb{T}^0 = 0$.

	Then the functor $\mathrm{MC}(\mb{T},-): \mc{C}^\dagger \to \mathrm{Set}$,
	$R \mapsto \mathrm{MC}(\mb{T},R)$
	is corepresentable by an object of $\mc{C}^\dagger$, $R^{\mathrm{univ}}$.
	Furthermore, $R^{\mathrm{univ}}$ can be written as a quotient of $\mb{Q}_p\{X_1,\ldots, X_g\}$,
    where $g := \dim_{\mb{Q}_p}\left(H^1(\mb{T})\right)$, by an ideal generated by $\dim_{\mb{Q}_p}\left(H^2(\mb{T})\right)$ elements.
\end{theorem}

\begin{remark}
	This is \cite[Theorem 1.4]{LiOcdt} with the following simplifications applied:
	\begin{enumerate}
		\item Typically, one can only represent $\mrm{MC}$ \emph{modulo gauge equivalence}
		as $\mrm{MC}$ itself is not a homotopy invariant.
		However, as $\mb{T}^0 = 0$, the gauge equivalence relation is the trivial equivalence relation.
		\item Since $\mb{T}$ is termwise Banach,
		the existence of a bounded retract follows from the finite-dimensionality of $H^\bullet(\mb{T})$
		by the open mapping theorem (\cite[Remark 2.28]{LiOcdt}).
%		\item It is crucial we are defining $H^i(\mb{T})$ as $Z^i(\mb{T})/B^i(\mb{T})$ \emph{without} taking closure of the boundaries.
%		In general this theorem would require a hypothesis on the existence of a ``bounded retract''
%		but for Banach spaces with finite-dimensional cohomology existence of such a retract follows from the Open Mapping Theorem
%		(\cite[Remark 2.19]{LiOcdt}).
%		\item The crux of the proof is making sense of explicit formulas for homotopy transfer in a rigid analytic setting.
%		In the Banach case the bijection between $\Hom(R^{\mathrm{univ}}, R)$
%		and $\mathrm{MC}(\mb{T},R)$ is in terms of these and hence somewhat explicit.
	\end{enumerate}
\end{remark}

To relate to the introduction of this section, when $H^2(\mathbb{T})=0$, $R^{\mrm{univ}}$ is the overconvergent symmetric algebra on $(H^1(\mathbb{T}))^{\vee}$
and
$$
\Hom(R^{\mrm{univ}}, A) \approx
\MC((H^{\bullet}(\mathbb{T})^\vee)^{\vee}, A)
\cong
\MC(H^{\bullet}(\mathbb{T}), A)
\cong
\MC(\mb{T}, A)
$$
(though this should be taken informally as $\MC(H^{\bullet}(\mathbb{T}), A)$ has not been defined).

\subsection{Algebraic incarnation}\label{ss.alg-inc}

We now prove an algebraic version of Corollary \ref{Cor.Germ}.
Recall the definition
$\Theta: (B^+_\crys)^{\varphi^n=p} \to C^n$,
$\Theta(b) := (\theta(b), \theta(\varphi(b)), \ldots, \theta(\varphi^{n-1}(b)))$.
Also recall the well-known fact that $\BCfrac{1}{n} \cong (B^+_\crys)^{\varphi^n=p}$
and, with a suitable choice of basis, $\iota = \Theta$.

For consistency, we will henceforth use $\BCfrac{1}{n}$
and $\iota$.
Evaluations $F \circ \iota$ are shorthand for
$\BCfrac{1}{n} \cong (B^+_\crys)^{\varphi^n=p} \xrightarrow{\Theta}
C^n \xrightarrow{F} C$.

In the following, $\iota$ has been extended $\mrm{Nil}^\dagger(R)$-linearly.
%In the following, $\theta$, $\varphi$ and $\Theta$
%are all extended $R$-linearly to
%$\mc{O}\left(\mc{Y}_{\mathrm{FF}}\right)^{\varphi^n=p} \otimes \mathrm{Nil}^\dagger(R)$.

\begin{proposition}
	\label[proposition]{prop:AlgebraicPadicSunscreen}
	Let $n$ and $F$ be as in Corollary \ref{Cor.Germ}.
	Then, for $R \in \mc{C}^\dg$,
	there is a natural (in $R$) isomorphism
	\begin{align*}
	\Hom(\mb{Q}_p\{X_1,\ldots,X_n\}, R)
	\cong
	\left\{x \in
    %(B^+_\crys)^{\varphi^n=p}
    %\mc{O}\left(\mc{Y}_{\mathrm{FF}}\right)^{\varphi^n=p}
    \BCfrac{1}{n}
		\otimes \mathrm{Nil}^{\dg}(R)
	\mid
	F\left(
        \iota(x)
		%\theta(x), \ldots, \theta(\varphi^{n-1}(x))
	\right) = 0
	\right\}.
	\end{align*}
\end{proposition}
\begin{proof}
	The proof is cleanest when $F$ is quadratic so we first assume this.
	Split $F = F_1 + F_2$ into its linear and quadratic part.

	Define $\mb{T}$ to be the DGA concentrated in degrees $[1,2]$:
	$$
	\mb{T} :=
	\left[
        \BCfrac{1}{n}
		%\mc{O}\left(
		%	\mc{Y}_{\mathrm{FF}}
		%\right)^{\varphi^n = p}
		\xrightarrow{F_1 \circ \iota}
		C
	\right],
	$$
	with multiplication
	$$
	x \cdot y :=
	-F_2'\left(
		\iota(x),
		\iota(y)
	\right),
	$$
	where $F_2'$ is any bilinear form such that $F_2' \circ \Delta = F_2$
	(and $\Delta: C^n \to (C^n)^{\otimes 2}$ denotes the diagonal map $v \mapsto v \otimes v$).

	Then, by construction,
	\begin{align*}
	\MC(\mb{T}, R) &=
	\{
		x \in \BCfrac{1}{n}
		\otimes \mrm{Nil}^\dagger(R)
		\mid
		F_1(\iota(x)) + F_2'(\iota(x),\iota(x)) = 0
	\} \\
	&=
	\{
		x \in \BCfrac{1}{n}
		\otimes \mrm{Nil}^\dagger(R)
		\mid
        F(\iota(x)) = 0
	\}.
	\end{align*}
	Hence the desired statement follows
	from \Cref{thm:Representability}.
	The conditions of the theorem are satisfied because:
	\begin{enumerate}
		\item $H^\bullet(\mb{T}) \cong [\mb{Q}_p^n \to 0]$ (in degrees 1 and 2) by Example \ref{example.bc1n}, and
		\item $\BCfrac{1}{n}$ has the structure of a Banach space as in \S\ref{sss.bc-banach}.
	\end{enumerate}

	The quadratic case is easier because \cite{LiOcdt}
	assumes $\mb{T}$ is a DGA.
	For the general case there are two options:
	the more conceptual is to adapt \cite{LiOcdt} to allow for $\mb{T}$
	to be a \emph{Fr{\'e}chet $A_\infty$-algebra}.
	We sketch this in the rest of this proof.
	Alternatively, one expects from the classical setting (\cite[Theorem 11.4.7]{LodayValletteAlgebraicOperads})
	that any $A_\infty$-algebra should be rectifiable
	into a homotopy equivalent DGA.
	In this case, we can explicitly construct a DGA whose Maurer--Cartan equation is the RHS of the proposition.
	Since this may seem ad hoc, we have opted to
	place it in \S\ref{appendix:DGAConstruction}.

	Write $F = F_1 + F_2 + \ldots$, with $F_m$ homogeneous of degree $m$.
	Continue to let $\mb{T}$ be the above chain complex.
	Let $F'_m:
	\left(C^n\right)^{\otimes m}
	\to C$
	be the symmetric multilinear form such that
	$F'_m \circ \Delta_m = F_m$,
	where $\Delta_m: C^n \to (C^n)^{\otimes m}$ is the diagonal morphism.
	The fact that $F'_m$ is symmetric is unimportant;
	any form such that $F'_m \circ \Delta_m = F_m$
	and that does not have unreasonably large denominators will work.

	Let $A^{\bullet} := \mrm{Free}^{\Fre}(\mb{T}[1])$ (\cite[Definition 2.4]{LiOcdt},
    except just as a graded module ignoring the differential).
	This is essentially the submodule of $\prod_{i=1}^{\infty} (\mb{T}[1])^{\otimes i}$
	of sequences whose valuations grow faster than linearly as $i \to \infty$
	(see \cite[Lemma 2.7]{LiOcdt} for a similar characterization).
	Let $m_i := F'_i \circ \iota^{\otimes i}$.
	On $A^0$ we can define $d' := \prod_{i=1}^{\infty} m_i$,
	meaning $d' = F'_i \circ \iota^{\otimes i}$ when restricted to $(\mb{T}[1])^{\otimes i}$
	and we take the infinite sum of these
	(this sum converges due to the growth rate condition on $\mrm{Free}^{\Fre}$
	and the assumption that $F$ converges in some neighborhood of $0$).
	On $A^j$ for $j \neq 0$, we define $d'$ to be 0.
	Define $d_2$ to be the coderivation of $A^{\bullet}$ extending $d'$ (\cite[Definition 2.17]{LiOcdt}).
	Then $(A^{\bullet}, d_2)$ is an $A_\infty$-algebra.
	
	The homotopy transfer of \textit{loc. cit.} can be easily adapted to apply to $A^{\bullet}$
	to construct the $A_\infty$-algebra denoted
	$H_\infty$ at the end of \cite[\S2.6]{LiOcdt}.
	By construction, $\mb{Q}_p \oplus (H_\infty)^{0,\vee,\mrm{ab}} \cong \mb{Q}_p\{X_1,\ldots,X_n\}$,
	and by arguing as in \cite[\S2.9]{LiOcdt}
	we see that
	$$
	\Hom(\mb{Q}_p \oplus (H_\infty)^{0,\vee,\mrm{ab}}, R)
	\cong
	\mrm{MC}\left(
		(A^{\bullet}, d_2),
		R
	\right),
	$$
	where
	\begin{equation}
    \label{eq:MCAinfAlgebra}
	\mrm{MC}\left(
		(A^{\bullet}, d_2),
		R
	\right)
	=
	\left\{
		x \in \mb{T}^1 \otimes \mrm{Nil}^\dagger(R)
		\bigm\vert
		\sum_{i=1}^{\infty} (m_i \otimes \mu_i)(x^{\otimes i}) = 0
	\right\}
	\end{equation}
	($\mu_i: R^{\otimes i} \to R$ denotes multiplication).
	Unraveling definitions, this is the RHS of the proposition.
\end{proof}

\subsection{Proof of the main result}\label{ss.top-version}

In this subsection, for concreteness, we let $\mb{T}$ be the DGA constructed in \S\ref{appendix:DGAConstruction}.
%(we note that one can modify the arguments to work for the $A_\infty$-algebra in the proof of \Cref{prop:AlgebraicPadicSunscreen}).

In an abuse of notation,
let $\MC(\mb{T}, \mb{Q}_p) := \{x \in \mb{T}^1 \mid dx - x \cdot x = 0\}$,
considered as a topological space
with the subspace topology from $\mb{T}^1$. 
% Define \MC(\mb{T}, \mb{Q}_p) $$ ... $$ (with notation as in \cref{eq:MCAinfAlgebra}.

Let $R^{\mathrm{univ}} := \mb{Q}_p\{X_1,\ldots,X_n\}$
denote the (co-)representing ring of \Cref{prop:AlgebraicPadicSunscreen}
and let $A^{\mrm{univ}} \in \MC(\mb{T}, R^{\mathrm{univ}})$ be the universal solution
(i.e., the element corresponding to $\id \in$ $\mrm{End}(R^{\mathrm{univ}})$ under the isomorphism of \Cref{prop:AlgebraicPadicSunscreen}).

Let $R_k^{\mathrm{univ}} := \mb{Q}_p\langle X_1/p^k, \ldots, X_n/p^k\rangle$.
Then, as $\MC$ commutes with filtered colimits,
there exists an $N$ such that
$A^{\mrm{univ}} \in \MC(\mb{T}, R_N^{\mathrm{univ}})$. 
Fix such an $N$. % and we now let $F$ denote this element.

Define $A: \left(p^N \mb{Z}_p\right)^n \to \MC(\mb{T}, \mb{Q}_p)$
as follows:
for $\mathbf{x} := (\mathbf{x}_1,\ldots,\mathbf{x}_n) \in \left(p^N \mb{Z}_p\right)^n$, let $A(\mathbf{x})$
be the image of $A^{\mrm{univ}}$ under the map
$\MC(\mb{T}, R_N^{\mathrm{univ}})
\xrightarrow{X_i \mapsto \mathbf{x}_i}
\MC(\mb{T}, \mb{Q}_p)$
base changed from $R_N^{\mathrm{univ}} \to \mb{Q}_p$,
sending $X_i$ to $\mathbf{x}_i$.

\begin{proposition}
	The morphism $A$
	induces a homeomorphism between
	some neighborhood of 0 in $\left(p^N \mb{Z}_p\right)^n$
	and some neighborhood of 0 in $\MC(\mb{T}, \mb{Q}_p)$.

	Hence, there exists a neighborhood of $0$
	in $\{x \in \BCfrac{1}{n} \mid F(\iota(x)) = 0\}$
	that is homeomorphic to $\mathbb{Z}_p^n$.
\end{proposition}
\begin{proof}
	We use a bit more about the construction behind \Cref{thm:Representability}.
	The construction involves choosing a continuous right inverse of $d^1$ (on its image), say $r^2: \mb{T}^2 \to \mb{T}^1$
	(which exists by the open mapping theorem and Lemma \ref{lem.splitting}).
	This gives a splitting $\mb{T}^1 \cong Z^1 \oplus \im(r^2)$.
	Let $p_{H^1}$ denote projection onto the first factor followed by projecting to $H^1$.
	The choice of presentation $R^{\mathrm{univ}} \cong \mb{Q}_p\{ X_1, \ldots, X_n \}$
	arises from a choice of basis, say $\mathbf{e}_1, \ldots, \mathbf{e}_n$ of $H^1$.
    These choices also determine the isomorphism in \Cref{prop:AlgebraicPadicSunscreen}.
    In particular,
	$$
	(\id \otimes p_{H^1})(A^{\mrm{univ}}) = X_1 \otimes \mathbf{e}_1 + \ldots + X_n \otimes \mathbf{e}_n \in R_N^{\mathrm{univ}} \otimes H^1.
	$$
	This follows from the construction of \Cref{thm:Representability} in \cite{LiOcdt}:
    the key fact is that, for $m > 1$ and the morphisms $\iota_{\infty,m}$ defined
    in \S2.7 of \textit{loc. cit.},
    we have $p_{H^1} \circ \iota_{\infty,m} = 0$
    as $\iota_{\infty,m}$ factors through $r^2$.
	
	Now, the key claim is that, under the identification
	$H^1 \cong (\mb{Q}_p)^n$ arising from the basis $\mathbf{e}_1, \ldots, \mathbf{e}_n$,
	$$
	A:
	H^1 \cong (\mb{Q}_p)^n \lrarrow \MC(\mb{T}, \mb{Q}_p):
	p_{H^1}
	$$
	are inverse homeomorphisms when restricted to sufficiently small neighborhoods of 0.
	By the above, $p_{H^1} \circ A = \id$.
	It remains to show that $A \circ p_{H^1} = \id$.

	Fix an $\alpha \in \MC(\mb{T}, \mb{Q}_p)$
    such that $p_{H^1}(\alpha) \in (p^N \mb{Z}_p)^n$.
	Let $\beta := (A \circ p_{H^1})(\alpha)$
	and $\varepsilon := \beta - \alpha$.
	Both $\alpha$ and $\beta$ satisfy the Maurer--Cartan equation:
	\begin{gather*}
	d\alpha - \alpha\cdot\alpha = 0,\\
	d(\alpha+\varepsilon) - (\alpha+\varepsilon)\cdot(\alpha+\varepsilon) = 0.
	\end{gather*}
	Subtract the two:
	$$
	d\varepsilon = (\alpha+\varepsilon)\cdot\varepsilon + \varepsilon\cdot\alpha.
	$$
	Observe that $\varepsilon \in \im(r^2)$, since $p_{H^1}(\alpha) = p_{H^1}(\beta)$.
	By the open mapping theorem, $d^1|_{\im(r^2)}: \im(r^2) \to \mb{T}^2$ is a topological isomorphism.
	Hence, $|d\varepsilon| \geq C|\varepsilon|$ for some $C$ independent of $\alpha$.
	By the definition of $\beta$, $|(\alpha+\varepsilon)| = |\beta| \leq C_2|\alpha|$ for some $C_2$ independent of $\alpha$.
	Hence,
	$$
	C |\varepsilon| \leq |d\varepsilon| =
	\lvert (\alpha+\varepsilon)\cdot\varepsilon + \varepsilon\cdot\alpha \rvert
	\leq C_\mu \max(1,C_2) |\alpha| |\varepsilon|,
	$$
    where $C_\mu$ is the norm of the multiplication map $\mb{T}^1 \otimes \mb{T}^1 \to \mb{T}^2$.
	We see that if $|\alpha| < C/\left(C_\mu \max(1,C_2)\right)$, then $\varepsilon = 0$, as desired.

	%The map $T$ is clearly continuous.
	%We now show the map is injective.
	%First, $\mathbf{x} \to A(\mathbf{x})$ is injective
	%as $(\proj_{Z^1} \otimes \id)$ is essentially a left inverse.
	%Next, $\proj_{\mrm{BC}(\mc{O}(1/n))}$ is injective on $\MC(\mb{T},\mb{Q}_p)$
	%by construction of $\mb{T}$.
	%Also, $\Theta$ is injective as any $x$ in the kernel vanishes at $\infty$,
	%hence gives a section of $\mc{O}(1/n)(\infty) \cong \mc{O}(-(n-1)/n)$ on $\mc{X}_{\mrm{FF}}$.
	%But $\Gamma(\mc{X}_{\mathrm{FF}}, \mc{O}(-(n-1)/n)) = 0$.
%
	%We now show surjectivity in a neighborhood.
	%Let $\alpha \in \im(\Theta) \cap \{F = 0\} \cap \mb{D}(0,1)$,
	%where $\mb{D}(0,r)$ is the open ball around 0 of radius $r$.

    The final claim follows
    since $\mb{T}$ is constructed so that
    $\MC(\mb{T}, \mb{Q}_p) \cong \{x \in \BCfrac{1}{n} \mid F(\iota(x)) = 0\}$
    in a neighborhood of $0$ where $F(\iota(x))$ converges.
\end{proof}

\subsection{DGA construction for general \texorpdfstring{$F$}{F}}
\label{appendix:DGAConstruction}

In this subsection we construct, for $F$ as in Corollary \ref{Cor.Germ} (that is, $F \in C[[x_1,\ldots,x_n]]$ with constant term $0$, non-vanishing gradient at $0$, and convergent in a neighborhood of $0$), a DGA $\mb{T}$ such that
$$
\MC(\mb{T}, R)
\cong
\{x \in \BCfrac{1}{n} \otimes \mrm{Nil}^\dagger(R) \mid F(\iota(x)) = 0\}.
$$
%We do the algebraic step more generally now.
%Like in the $x-y^2=0$ example the linear terms are encoded in the differential
%and the multiplication encodes the degree 2 terms.
%The trick with the higher degree terms is to turn them into iterated multiplications using extra variables.

Again, $\mb{T}$ is in degrees $[1,2]$.
Fix an $r \in \mb{Z}$ such that $F$ converges on
$\{|X_1|,\ldots, |X_n| \leq p^{-r}\}$.
Fix indeterminate variables ``$X_1, \ldots, X_n, Y$''
and define
\begin{gather*}
	\mc{P} := \widehat{\bigoplus_{\substack{a_1,\ldots,a_n \geq 0 \\ a_1 + \ldots + a_n > 0}}} \mc{O}_{C}\cdot (p^rX_1)^{a_1}\cdots (p^rX_n)^{a_n}\left[\frac{1}{p}\right], \\
	\mb{T}^1 := \BCfrac{1}{n} \oplus \mc{P}, \\
    \mb{T}^2 := \mc{P}
    \oplus C \cdot Y.
\end{gather*}

The differential is as follows.
For $a \in \BCfrac{1}{n}$,
$$
d^1(a)
:=
-\bigoplus_{i=1}^n \theta(\varphi^{i-1}(a)) \cdot X_i.
$$
On $\mc{P}$,
$d^1|_{\mc{P}} := \id_{\mc{P}} \oplus G$,
where $G(X_1^{a_1}\cdots X_n^{a_n})$
is the coefficient of
$X_1^{a_1}\cdots X_n^{a_n}$ in $F$.

Clearly $\proj_{\mc{P}} \circ d^1$ is surjective
with kernel $\cong \BCfrac{1}{n}$.
By Example \ref{example.bc1n} and the assumption
that $F$ has non-vanishing gradient,
it is not hard to see
that $d^1$ is surjective and $H^1 \cong \mb{Q}_p^n$.

\smallskip
We now define the multiplication, $\mu: \mb{T}^1 \otimes \mb{T}^1 \to \mb{T}^2$.
First,
\begin{align*}
    &\proj_{Y} \circ \mu = 0, \\
    &\proj_{X_1^{a_1}\cdots X_n^{a_n}} \circ \mu = 0,
    \quad \text{if $a_1+\ldots+a_n = 1$}.
\end{align*}
For $a_1+\ldots+a_n > 1$,
let $i$ be the first non-zero $a_i$.
Then
\begin{align*}
\proj_{X_1^{a_1}\cdots X_n^{a_n}} &\circ \mu := \\
&\mb{T}^1 \otimes \mb{T}^1 \twoheadrightarrow
\BCfrac{1}{n} \otimes
\left(
    C \cdot X_1^{a_1}\cdots X_i^{a_i-1} \cdots X_n^{a_n}
\right)
\to C \cdot X_1^{a_1}\cdots X_n^{a_n},
\end{align*}
where the last arrow is
$$
x \otimes
\left(
    c\cdot X_1^{a_1}\cdots X_i^{a_i-1} \cdots X_n^{a_n}
\right)
\mapsto \theta(\varphi^{i-1}(x))c \cdot X_1^{a_1}\cdots X_n^{a_n}.
$$

\smallskip

One can now check that $\proj_{\mrm{BC}(1/n)}$
induces an isomorphism
$$
\MC(\mb{T}, R)
\cong
\{x \in \BCfrac{1}{n} \otimes \mrm{Nil}^\dagger(R) \mid F(\iota(x)) = 0\}.
$$
In essence, any $x \in \MC(\mb{T}, R)$
is determined by its projection to $\BCfrac{1}{n}$
because $d^1$ acts as the identity on the $\mc{P}$ component
and $\mu$ raises degrees.

\section{A non-linear example}\label{s.non-linear}
In this section we explain how to equip the $C$-points of the infinite level Lubin--Tate space with a natural structure of a $\mathbb{Q}_p$-analytic Banach manifold (\S\ref{s.ilt-banach}). We then explain that the Hodge and Hodge--Tate period maps are analytic; the tangent spaces admit natural structures as Banach--Colmez spaces, and the derivatives in this theory match previous computations in the infinitesimal theories of \cite{Howe.CohomologicalSmoothnessConjecture, IvanovWeinstein.TheSmoothLocusInInfiniteLevelRapoportZinkSpaces}. In particular, we can then take further transverse intersections with analytic submanifolds of a product of projective spaces, enriching the geometry of the transverse intersections studied in \cite[\S6]{HoweKlevdal.AdmissiblePairsAndpAdicHodgeStructuresIIIVariationAndUnlikelyIntersection} in the context of $p$-adic unlikely intersections (see Remark \ref{remark.hk-intersections}). 

\subsection{The infinite level Lubin--Tate space as a Banach manifold}\label{s.ilt-banach}
For any $n > 1$, there is a natural identification $\wedge^n \mathcal{O}(\frac{1}{n})=\mathcal{O}(1)$, inducing a map
\[ \BCfrac{1}{n}^n \xrightarrow{\det_\BC} \BC(1). \]
Let $\mathcal{P}_C: \BCfrac{1}{n}^n \rightarrow M_{n \times n}(C)$ denote the map obtained by specializing at $\infty$ (more precisely, by viewing an element of $\BCfrac{1}{n}^n$ as a homomorphism from $\mathcal{O}^n$ to $\mathcal{O}(\frac{1}{n})$ and then using the canonical trivializations over $\infty$). Let $U \subseteq M_{n\times n}(C)$ be the open subset of matrices of rank at least $n-1$.

\begin{lemma} 
The map $\det_\BC$ is an analytic map of $\mathbb{Q}_p$-analytic Banach manifolds, and it is a submersion (i.e., its derivative is surjective on tangent spaces) after restriction to the open submanifold $V:=\mathcal{P}_C^{-1}(U)$. 
\end{lemma}
\begin{proof}
Let $\mathcal{P}_{B^+_\dR/t^2}: \BCfrac{1}{n}^n \rightarrow M_{n \times n}(B^+_\dR/t^2)$ denote the map obtained by specializing to the first infinitesimal neighborhood of $\infty$. We obtain a commutative diagram
% https://q.uiver.app/#q=WzAsNCxbMCwwLCJcXEJDKDEvbilebiJdLFswLDEsIk1fe24gXFx0aW1lcyBufShCXitfXFxkUi90XjIpIl0sWzEsMCwiXFxCQygxKSJdLFsxLDEsIkJeK19cXGRSL3ReMiJdLFswLDIsIlxcZGV0Il0sWzAsMV0sWzIsM10sWzEsMywiXFxkZXQiLDJdXQ==
\[\begin{tikzcd}
	{\BCfrac{1}{n}^n} & {\BC(1)} \\
	{M_{n \times n}(B^+_\dR/t^2)} & {B^+_\dR/t^2}
	\arrow["\det_\BC", from=1-1, to=1-2]
	\arrow["\mathcal{P}_{B^+_\dR/t^2}", from=1-1, to=2-1]
	\arrow[from=1-2, to=2-2]
	\arrow["\det"', from=2-1, to=2-2]
\end{tikzcd}\]
where the bottom horizontal arrow is matrix determinant and the vertical arrows are both the inclusions of closed subspaces of $\mathbb{Q}_p$-Banach spaces. It follows that $\det_\BC$ is analytic --- indeed, $\det \circ \mathcal{P}_{B^+_\dR/t^2}$ is a homogeneous polynomial function of degree $n$, thus so is $\det_\BC$ as it is recovered by composing with a projection onto the closed subspace $\BC(1)$ of $B^+_\dR/t^2$ (which exists as any closed subspace is complemented by \cite[Proposition 10.5]{Schneider.NonarchimedeanFunctionalAnalysis}).

We now establish the submersion claim. To that end, suppose $v=(v_1, \ldots, v_n) \in V$, and write $\overline{v}_i$ for the image of $v_i$ in $C^n$. Then, up to permutation, we can assume  $\overline{v}_1, \ldots, \overline{v}_{n-1}$  are linearly independent in $C^n$ (the $\overline{v}_i$ are the columns of $\mathcal{P}_C(v)$). We consider $\BCfrac{1}{n}$ as a subspace of $\BCfrac{1}{n}^n=T_{V,v}$ by $w \mapsto (0, \ldots, 0, w)$. Applying $d\det_{\BC,v}$ to $(0, \ldots, 0, w)$ and then projecting from $T_{\BC(1), \det_\BC(v)}=\BC(1)$ to $C$ yields $\det(\overline{v}_1, \overline{v}_2, \ldots, \overline{v}_{n-1},\overline{w})$. Since $\BCfrac{1}{n}$ surjects onto $C^n/(C\overline{v}_1 + \ldots + C\overline{v}_{n-1})$ by Example \ref{example.bc1n}, we find the map $\BCfrac{1}{n} \rightarrow \BC(1)$ is a surjection after composition with $\BC(1) \rightarrow C$. We claim this implies $\BCfrac{1}{n} \rightarrow \BC(1)$ is a surjection: indeed, the map in question comes from a map of vector bundles $\mathcal{O}(\frac{1}{n}) \rightarrow \mathcal{O}(1)$ that fits into a map of modifications

\[\begin{tikzcd}
	0 & {\mathcal{O}^n} & {\mathcal{O}(\frac{1}{n})} & {\infty_*(C^n/(C\overline{v}_1+\ldots +C\overline{v}_{n-1}))} & 0 \\
	0 & {\mathcal{O}} & {\mathcal{O}(1)} & {\infty_*C} & 0
	\arrow[from=1-1, to=1-2]
	\arrow[from=1-2, to=1-3]
	\arrow[from=1-2, to=2-2]
	\arrow[from=1-3, to=1-4]
	\arrow[from=1-3, to=2-3]
	\arrow[from=1-4, to=1-5]
	\arrow[from=1-4, to=2-4]
	\arrow[from=2-1, to=2-2]
	\arrow[from=2-2, to=2-3]
	\arrow[from=2-3, to=2-4]
	\arrow[from=2-4, to=2-5]
\end{tikzcd}\]
and the right vertical arrow is surjective. If the left vertical arrow were zero, we would obtain a splitting of the bottom exact sequence, which does not exist; thus the left vertical arrow is non-zero so surjective. Both sequences remain exact after taking global sections, and since the left and right vertical arrows are surjective on global sections, so is the middle one. 
\end{proof}

Letting $\mathscr{M}_\infty$ denote the height $n$ infinite level Lubin--Tate space, it follows from the results of \cite{ScholzeWeinstein.ModuliOfpDivisibleGroups} that $\mathscr{M}_\infty(C)$ is naturally identified with the pre-image under $\det_\BC$ of 
\[ \mathbb{Q}_p(1)\backslash\{0\}=\mathbb{Q}_p^\times \cdot t \subseteq \BC(1).\]

Note that $\det_{\BC}^{-1}(\BC(1) \backslash \{0\}) \subseteq V$ --- indeed, if the $C$-period matrix of $v \in \BCfrac{1}{n}^n$ has rank $< n - 1$ it forces $\det_\BC(v)$ to be divisible by $t^2$ and thus equal to $0$ as $\BC(1) \hookrightarrow B^+_\dR/t^2$. By the usual argument, the implicit function theorem (Theorem \ref{thm.implicit-function}) implies that the pre-image of a submanifold under a submersion is a submanifold (again using that any closed subspace of a $\mathbb{Q}_p$-Banach space is complemented \cite[Proposition 10.5]{Schneider.NonarchimedeanFunctionalAnalysis}), thus we find that $\mathscr{M}_\infty(C)$ has a natural $\mathbb{Q}_p$-analytic Banach manifold structure as a closed submanifold of $\det_\BC^{-1}(\BC(1) \backslash\{0\})$ (which itself is open in $\BCfrac{1}{n}^n$).

\subsection{Tangent spaces and derivatives of period maps}
As an immediate consequence of the construction, the $\mathbb{Q}_p$-Banach tangent space of $\mathscr{M}_\infty(C)$ at any point has a natural structure of a Banach--Colmez space (associated to a modification of the vector bundle $\mathcal{O}(\frac{1}{n})^{\oplus n}$), and these are readily seen to match with the Banach--Colmez Tangent Spaces coming from nilpotent structures either via the Fargues--Scholze Jacobian criterion as in \cite{IvanovWeinstein.TheSmoothLocusInInfiniteLevelRapoportZinkSpaces} or the theory of inscribed $v$-sheaves as in \cite{Howe.CohomologicalSmoothnessConjecture} (in the former case one should in fact take a connected component of $\mathscr{M}_\infty$ corresponding to the fiber over a single point in $\mathbb{Q}_p(1)\backslash\{0\}$; the inscribed computation makes sense for both variants). The advantage of the present computation is that it also establishes a link between the underlying $\mathbb{Q}_p$-Banach space and the local geometry of the $C$-points --- in a small enough open neighborhood, they are the same. 

One of the key points in the inscribed computation \cite{Howe.CohomologicalSmoothnessConjecture} is that the Hodge and Hodge--Tate period maps are differentiable. In the present setting, they are both analytic maps of $\mathbb{Q}_p$-analytic Banach manifolds: indeed, the Hodge--Tate/Hodge period map is the composition of the period matrix $\mathcal{P}_C$ with the map sending an $n \times n$ matrix of rank $n-1$ to its kernel/cokernel viewed as a point in $\mathbb{P}^{n-1}(C)$. It is readily verified that the differentials of these maps also match the differentials computed in the inscribed setting. 

\begin{remark}\label{remark.hk-intersections}
Consider the variant $\mathscr{M}_\infty^\circ(C)$ obtained by taking the fiber over a fixed element of $\mathbb{Q}_p(1)\backslash\{0\}$. Transverse intersections of   $\mathscr{M}_\infty^\circ$ with smooth rigid analytic subvarieties of $\mathbb{P}^{n-1} \times \mathbb{P}^{n-1}$ are studied in \cite[\S6]{HoweKlevdal.AdmissiblePairsAndpAdicHodgeStructuresIIIVariationAndUnlikelyIntersection} via a formulation based on the theory of inscribed $v$-sheaves (here we mean relative transversality for the product of the two period maps). The nilpotent computations with inscribed $v$-sheaves all match the Banach theory so, at the level of $C$-points, these transverse intersections are $\mathbb{Q}_p$-analytic Banach manifolds that look locally like their inscribed tangent spaces. Thus the differential computations going into the description of exceptional/non-exceptional intersections are directly related to the punctual geometry of those intersections. We plan to return to this in a more general setting in future work. 
\end{remark}

\begin{remark}\label{remark.BanachManifolds}
    As in Question \ref{q.extension}, it is interesting to also consider how these Banach manifold structures vary with $C$. In this context, it is most natural to view $\mathscr{M}_\infty$ as a Banach Manifold, by which we mean roughly a $v$-sheaf on strictly totally disconnected perfectoid spaces equipped with a factorization through $\mathbb{Q}_p$-analytic Banach manifolds. This category, which we expect to give a common refinement of the theory of inscribed $v$-sheaves and the punctual $\mathbb{Q}_p$-analytic Banach manifold perspective of the present paper, is being developed elsewhere. 
\end{remark}

\bibliographystyle{plain}
\bibliography{references, preprints, ray}

@article{Howe.CohomologicalSmoothnessConjecture,
title={A cohomological smoothness conjecture for moduli of mixed characteristic local shtukas with one leg},
author={Howe, Sean},
journal={ar{X}iv:2508.11595},
year={2025}
}

@article{HoweKlevdal.AdmissiblePairsAndpAdicHodgeStructuresIIIVariationAndUnlikelyIntersection,
  author = {Howe, Sean and Klevdal, Christian},
  title  = {Admissible pairs and $p$-adic {H}odge structures {III}: Variation and unlikely intersection},
journal ={arXiv:2603.22610},
year={2026}
}

@article{Howe.InscriptionTwistorsAndPAdicPeriods,
AUTHOR={Howe, Sean},
TITLE={Inscription, twistors, and $p$-adic periods},
YEAR={2025},
JOURNAL={ar{X}iv:2508.11589}
}

@article{Howe.TheRelativisticPAdicSunscreenConjecture,
AUTHOR={Howe, Sean},
TITLE={The relativistic $p$-adic sunscreen conjecture},
YEAR={2026},
JOURNAL={ar{X}iv:2604.01245}
}

@book{BoschLectures,
    AUTHOR = {Bosch, Siegfried},
     TITLE = {Lectures on formal and rigid geometry},
    SERIES = {Lecture Notes in Mathematics},
    VOLUME = {2105},
 PUBLISHER = {Springer, Cham},
      YEAR = {2014},
     PAGES = {viii+254},
      ISBN = {978-3-319-04416-3; 978-3-319-04417-0},
   MRCLASS = {14G22 (11S85 14-02 14D15)},
  MRNUMBER = {3309387},
MRREVIEWER = {Alessandra\ Bertapelle},
       DOI = {10.1007/978-3-319-04417-0},
       URL = {https://doi.org/10.1007/978-3-319-04417-0},
}

@misc{LiOcdt,
    title={Overconvergent deformation theory and an application to $(\varphi,\Gamma)$-modules},
    author={Ray Li},
    howpublished={In preparation}
}

@book{LodayValletteAlgebraicOperads,
    AUTHOR = {Loday, Jean-Louis and Vallette, Bruno},
     TITLE = {Algebraic operads},
    SERIES = {Grundlehren der mathematischen Wissenschaften [Fundamental
              Principles of Mathematical Sciences]},
    VOLUME = {346},
 PUBLISHER = {Springer, Heidelberg},
      YEAR = {2012},
     PAGES = {xxiv+634},
      ISBN = {978-3-642-30361-6},
   MRCLASS = {18D50 (16E99)},
  MRNUMBER = {2954392},
MRREVIEWER = {Andrey\ Yu.\ Lazarev},
       DOI = {10.1007/978-3-642-30362-3},
       URL = {https://doi.org/10.1007/978-3-642-30362-3},
}

@book{Schneider.NonarchimedeanFunctionalAnalysis,
    AUTHOR = {Schneider, Peter},
     TITLE = {Nonarchimedean functional analysis},
    SERIES = {Springer Monographs in Mathematics},
 PUBLISHER = {Springer-Verlag, Berlin},
      YEAR = {2002},
     PAGES = {vi+156},
      ISBN = {3-540-42533-0},
   MRCLASS = {46S10 (46-01 47S10)},
  MRNUMBER = {1869547},
MRREVIEWER = {Dinam\'erico\ P.\ Pombo, Jr.},
       DOI = {10.1007/978-3-662-04728-6},
       URL = {https://doi-org.ezproxy.lib.utah.edu/10.1007/978-3-662-04728-6},
}

@book {Serre.LieAlgebrasAndLieGroups,
    AUTHOR = {Serre, Jean-Pierre},
     TITLE = {Lie algebras and {L}ie groups},
    SERIES = {Lecture Notes in Mathematics},
    VOLUME = {1500},
      NOTE = {1964 lectures given at Harvard University,
              Corrected fifth printing of the second (1992) edition},
 PUBLISHER = {Springer-Verlag, Berlin},
      YEAR = {2006},
     PAGES = {viii+168},
      ISBN = {978-3-540-55008-2; 3-540-55008-9},
   MRCLASS = {17-01 (22-01)},
  MRNUMBER = {2179691},
}

@Article{ScholzeWeinstein.ModuliOfpDivisibleGroups,
  author     = {Scholze, Peter and Weinstein, Jared},
  journal    = {Camb. J. Math.},
  title      = {Moduli of {$p$}-divisible groups},
  year       = {2013},
  issn       = {2168-0930},
  number     = {2},
  pages      = {145--237},
  volume     = {1},
  doi        = {10.4310/CJM.2013.v1.n2.a1},
  fjournal   = {Cambridge Journal of Mathematics},
  mrclass    = {14L05 (11G25 14C30 14D20)},
  mrnumber   = {3272049},
  mrreviewer = {Rui Miguel Saramago},
  url        = {http://dx.doi.org/10.4310/CJM.2013.v1.n2.a1},
}

@Article{IvanovWeinstein.TheSmoothLocusInInfiniteLevelRapoportZinkSpaces,
  author     = {Ivanov, Alexander B. and Weinstein, Jared},
  journal    = {Compos. Math.},
  title      = {The smooth locus in infinite-level {R}apoport-{Z}ink spaces},
  year       = {2020},
  issn       = {0010-437X},
  number     = {9},
  pages      = {1846--1872},
  volume     = {156},
  doi        = {10.1112/s0010437x20007332},
  fjournal   = {Compositio Mathematica},
  mrclass    = {14G45 (14G22 14L05)},
  mrnumber   = {4170574},
  mrreviewer = {Hui June Zhu},
  url        = {https://doi.org/10.1112/s0010437x20007332},
}

@book {Bourbaki.VarietesDifferentiellesEtAnalytique1-7,
    AUTHOR = {Bourbaki, N.},
     TITLE = {\'El\'ements de math\'ematique. {F}asc. {XXXIII}.
              {V}ari\'et\'es diff\'erentielles et analytiques. {F}ascicule
              de r\'esultats ({P}aragraphes 1 \`a{} 7)},
    SERIES = {Actualit\'es Scientifiques et Industrielles [Current
              Scientific and Industrial Topics]},
    VOLUME = {No. 1333},
 PUBLISHER = {Hermann, Paris},
      YEAR = {1967},
     PAGES = {97},
   MRCLASS = {57.10 (14.00)},
  MRNUMBER = {219078},
}

\end{document}